\documentclass{article}%
\usepackage{amssymb}
\usepackage{amsmath}%
\usepackage{amsfonts}%
\usepackage{graphicx}
\providecommand{\U}[1]{\protect\rule{.1in}{.1in}}
\newtheorem{theorem}{Theorem}

\newtheorem{definition}[theorem]{Definition}

\newtheorem{lemma}[theorem]{Lemma}

\newtheorem{remark}[theorem]{Remark}

\begin{document}

\title{Existence of solitons in the nonlinear beam equation.}
\author{Vieri Benci$^{\ast}$, Donato Fortunato$^{\ast\ast}$\\$^{\ast}$Dipartimento di Matematica Applicata \textquotedblleft U.
Dini\textquotedblright\\Universit\`{a} degli Studi di Pisa \\Largo Bruno Pontecorvo 1/c, 56127 Pisa, Italy\\e-mail: benci@dma.unipi.it\\College of Science, Department of Mathematics \\King Saud University \\Riyadh, 11451, Saudi Arabia \\$^{\ast\ast}$Dipartimento di Matematica \\Universit\`{a} degli Studi di Bari Aldo Moro\\Via Orabona 4, 70125 Bari, Italy\\e-mail: fortunat@dm.uniba.it}
\maketitle

\begin{abstract}
This paper concerns with the existence of solitons, namely stable solitary
waves in the nonlinear beam equation (NBE) with a suitable nonlinearity. An
equation of this type has been introduced in \cite{mkw87} as a model of a
suspension bridge. We prove both the existence of solitary waves for a large
class of nonlinearities and their stability. As far as we know this is the
first result about stability of solitary waves in NBE.

\end{abstract}
\tableofcontents

\bigskip


\bigskip

\bigskip

\section{Introduction}

Let us consider the nonlinear beam equation%
\begin{equation}
\frac{\partial^{2}u}{\partial t^{2}}+\frac{\partial^{4}u}{\partial x^{4}%
}+W^{\prime}(u)=0\label{1}%
\end{equation}
where $u=u(t,x),\ $and $W\in C^{1}(\mathbb{R)}$. In this paper we will prove
that, under suitable assumptions, equation (\ref{1}) admits soliton solutions.
Roughly speaking a \textit{solitary wave} is a solution of a field equation
whose energy travels as a localized packet and which preserves this
localization in time. A \textit{soliton} is a solitary wave which exhibits
some form of stability so that it has a particle-like behavior (see e.g.
\cite{milano} or \cite{befogranas}). Following \cite{milano}, a soliton or
solitary wave is called \textit{hylomorphic} if its stability is due to a
particular ratio between \textit{energy }$E$ and the \textit{hylenic}
\textit{charge }$C$ which is another integral of motion. More precisely, a
soliton $\mathbf{u}_{0}$ is hylomorphic if%
\[
E(\mathbf{u}_{0})=\min\left\{  E(\mathbf{u})\ |\ C(\mathbf{u})=C(\mathbf{u}%
_{0})\right\}  .
\]
The physical meaning of $C$ depends on the problem (in this case $C$ is the
\textit{momentum,\ }see section \ref{smr}). The main result of this paper is
the proof of the existence of hylomorphic solitons for equation (\ref{1})
provided that $W$ satisfies suitable assumptions (namely (W-i), (W-ii) and
(W-iii) of section \ref{smr}). In particular, these assumptions are satisfied
by%
\begin{equation}
W(s)=\left\{
\begin{array}
[c]{cc}%
\frac{1}{2}s^{2} & for\ \ s\leq1\\
& \\
s-\frac{1}{2} & for\ \ s\geq1
\end{array}
\ \right. \label{W1}%
\end{equation}

Equation (\ref{1}) with $W(s)$ as in (\ref{W1}) has been proposed as model for
a suspension bridge (see \cite{mkw87}, \cite{lam}, \cite{lambis}). In
particular in \cite{McKW90} and \cite{SaWei09} the existence of travelling
waves has been proved.

Observe that $u(t,x)-1$ denotes the displacement of the beam from the unloaded
state $u(x)\equiv1$ and the bridge is seen as a vibrating beam supported by
cables which are treated as springs. The force relative to the potential
$W(s)$ in (\ref{W1}) is given by%
\[
F(s)=-W^{\prime}(s)=\left\{
\begin{array}
[c]{cc}%
-s & for\ \ s\leq1\\
& \\
-1 & for\ \ s\geq1;
\end{array}
\ ,\right.
\]
namely, for $s\geq1,$ only the costant gravity force $-1$ acts; while, for
$s\leq1,$ an elastic force (of intensity $1-s$), due to the suspension cables,
must be added to the costant gravity force $-1$ . Of course assumptions (W-i),
(W-ii) and (W-iii) are satisfied also by the potential
\begin{equation}
W(s)=s-1+e^{-s}\label{W2}%
\end{equation}
which has been considered in \cite{McKW90} and in \cite{SaWei09} as an
alternative smooth model for a suspension bridge.

\bigskip

\section{Hylomorphic solitary waves and solitons\label{HS}}

\subsection{An abstract definition of solitary waves and solitons\label{be}}

Solitary waves and solitons are particular \textit{states} of a dynamical
system described by one or more partial differential equations. Thus, we
assume that the states of this system are described by one or more
\textit{fields} which mathematically are represented by functions
\begin{equation}
\mathbf{u}:\mathbb{R}^{N}\rightarrow V\label{lilla}%
\end{equation}
where $V$ is a vector space with norm $\left\vert \ \cdot\ \right\vert _{V}$
which is called the internal parameters space. We assume the system to be
deterministic; this means that it can be described as a dynamical system
$\left(  X,\gamma\right)  $ where $X$ is the set of the states and
$\gamma:\mathbb{R}\times X\rightarrow X$ is the time evolution map. If
$\mathbf{u}_{0}(x)\in X,$ the evolution of the system will be described by the
function
\begin{equation}
\mathbf{u}\left(  t,x\right)  :=\gamma_{t}\mathbf{u}_{0}(x).\label{flusso}%
\end{equation}
We assume that the states of $X$ have "finite energy" so that they decay at
$\infty$ sufficiently fast.

We give a formal definition of solitary wave:

\begin{definition}
\label{solw} A state $\mathbf{u}(x)\in X$ is called solitary wave if there is
$\xi(t)$ such that
\[
\gamma_{t}\mathbf{u}(x)=\mathbf{u}(x-\xi(t)).
\]

\end{definition}

The solitons are solitary waves characterized by some form of stability. To
define them at this level of abstractness, we need to recall some well known
notions in the theory of dynamical systems.

\begin{definition}
A set $\Gamma\subset X$ is called \textit{invariant} if $\forall\mathbf{u}%
\in\Gamma,\forall t\in\mathbb{R},\ \gamma_{t}\mathbf{u}\in\Gamma.$
\end{definition}

\begin{definition}
Let $\left(  X,d\right)  $ be a metric space and let $\left(  X,\gamma\right)
$ be a dynamical system. An invariant set $\Gamma\subset X$ is called stable,
if $\forall\varepsilon>0,$ $\exists\delta>0,\;\forall\mathbf{u}\in X$,
\[
d(\mathbf{u},\Gamma)\leq\delta,
\]
implies that
\[
\forall t\geq0,\text{ }d(\gamma_{t}\mathbf{u,}\Gamma)\leq\varepsilon.
\]

\end{definition}

\bigskip

Let $G$ be the group induced by the translations in $\mathbb{R}^{N},$ namely,
for every $\tau\in\mathbb{R}^{N},\ $the transformation $g_{\tau}\in G$ is
defined as follows:
\begin{equation}
\left(  g_{\tau}\mathbf{u}\right)  \left(  x\right)  =\mathbf{u}\left(
x-\tau\right)  .\label{ggg}%
\end{equation}

\begin{definition}
A subset $\Gamma\subset X$ is called $G$-invariant if
\[
\forall\mathbf{u}\in\Gamma,\ \forall\tau\in\mathbb{R}^{N},\ g_{\tau}%
\mathbf{u}\in\Gamma.
\]

\end{definition}

\begin{definition}
A closed $G$-invariant set $\Gamma\subset X$ is called $G$-compact if for any
sequence $\mathbf{u}_{n}(x)$ in $\Gamma$ there is a sequence $\tau_{n}%
\in\mathbb{R}^{N},$ such that $\mathbf{u}_{n}(x-\tau_{n})$ has a converging subsequence.
\end{definition}

Now we are ready to give the definition of soliton:

\begin{definition}
\label{ds} A solitary wave $\mathbf{u}(x)$ is called soliton if there is an
invariant set $\Gamma$ such that

\begin{itemize}
\item (i) $\forall t,\ \gamma_{t}\mathbf{u}(x)\in\Gamma,$

\item (ii) $\Gamma$ is stable,

\item (iii) $\Gamma$ is $G$-compact.
\end{itemize}
\end{definition}

Usually, in the literature, the kind of stability described by the above
definition is called \textit{orbital stability}.

\begin{remark}
The above definition needs some explanation. For simplicity, we assume that
$\Gamma$ is a manifold (actually, this is the generic case in many
situations). Then (iii) implies that $\Gamma$ is finite dimensional. Since
$\Gamma$ is invariant, $\mathbf{u}_{0}\in\Gamma\Rightarrow\gamma_{t}%
\mathbf{u}_{0}\in\Gamma$ for every time. Thus, since $\Gamma$ is finite
dimensional, the evolution of $\mathbf{u}_{0}$ is described by a finite number
of parameters$.$ Thus the dynamical system $\left(  \Gamma,\gamma\right)
$\ behaves as a point in a finite dimensional phase space. By the stability of
$\Gamma$, a small perturbation of $\mathbf{u}_{0}$ remains close to $\Gamma.$
However, in this case, its evolution depends on an infinite number of
parameters. Thus, this system appears as a finite dimensional system with a
small perturbation. Since $\dim(G)=N$, $\dim\left(  \Gamma\right)  \geq N$ and
hence, the \textquotedblright state\textquotedblright\ of a soliton is
described by $N$ parameters which define its position and, may be, other
parameters which define its \textquotedblright internal
state\textquotedblright.
\end{remark}

\subsection{Integrals of motion and hylomorphic solitons\label{im}}

In recent papers (see e.g. \cite{milano}, \cite{hylo}, \cite{befo2011}), the
notion of \textit{hylomorphic soliton }has been introduced and analyzed. The
existence and the properties of hylomorphic solitons are guaranteed by the
interplay between the \textit{energy }$E$ and an other integral of motion
which, in the general case, is called \textit{hylenic charge} and it will be
denoted by\textit{\ }$C.$ More precisely:

\begin{definition}
\label{tdc}Assume that the dynamical system has two first integrals of motion
$E:X\rightarrow\mathbb{R}$ and $C:X\rightarrow\mathbb{R}$. A soliton
$\mathbf{u}_{0}\in X$ is hylomorphic if $\Gamma$ (as in Def. \ref{ds}) has the
following structure%
\[
\Gamma=\Gamma\left(  e_{0},p_{0}\right)  =\left\{  \mathbf{u}\in
X\ |\ E(\mathbf{u})=e_{0},\ C(\mathbf{u})=p_{0}\right\}
\]
where%
\[
e_{0}=\min\left\{  E(\mathbf{u})\ |\ C(\mathbf{u})=p_{0}\right\}
\]
for some $p_{0}\in\mathbb{R}$.
\end{definition}

Clearly, for a given $p_{0}$ the minimum of $E$ might not exist; moreover,
even if the minimum exists, it is possible that $\Gamma\ $does not satisfies
(ii) or (iii) of def. \ref{ds}.

In this section, we present an abstract theorem which guarantees the existence
of hylomorphic solitons. Before stating the abstract theorems, we need some definitions:

\begin{definition}
A functional $J$ on $X$ is called $G$\emph{-invariant} if
\[
\forall g\in G,\text{ }\forall\mathbf{u}\in X,\ J\left(  g\mathbf{u}\right)
=J\left(  \mathbf{u}\right)  .
\]

\end{definition}

\begin{definition}
Let $G$ be a group of tranlations acting on $X.$ A sequence $\mathbf{u}_{n}$
in $X$ is called $G$\emph{-compact }if we can extract a subsequence
$\mathbf{u}_{n_{k}}$ such that there exists a sequence $g_{k}\in G$ such that
$g_{k}\mathbf{u}_{n_{k}}$ is convergent.\ A functional $J$ on $X$ is called
$G$\emph{-compact} if any minimizing sequence of $J$ is\emph{\ }$G$-compact.
\end{definition}

\begin{remark}
Clearly, a $G$-compact functional admits a minimizer. Moreover, if $J$ is
$G$-invariant and $\mathbf{u}_{0}$ is a minimizers, then $\left\{
g\mathbf{u_{0}\ |\ }g\in G\right\}  $ is a set of minimizers; so, if $G$ is
not compact, the set of minimizers is not compact (unless $\mathbf{u_{0}}$ is
a constant). This fact adds an extra difficulty to this kind of problems.
\end{remark}

We make the following (abstract) assumptions on the dynamical system
$(X,\gamma)$:

\begin{itemize}
\item (EC-1) there are two first integrals $E:X\rightarrow\mathbb{R}$ and
$C:X\rightarrow\mathbb{R}.$

\item (EC-2) $E(\mathbf{u})$ and $C(\mathbf{u})$ are $G$-invariant.
\end{itemize}

\bigskip

\begin{theorem}
\label{astra} Assume that the dynamical system $(X,\gamma)$ satisfies (EC-1)
and (EC-2). Moreover we set%
\begin{equation}
J(\mathbf{u})=\frac{E(\mathbf{u})}{\left\vert C(\mathbf{u})\right\vert
}+\delta E(\mathbf{u})\label{jd}%
\end{equation}
where $\delta$ is a positive constant and assume that $J$ is $G$-compact. Then
$J(\mathbf{u})$ has a minimizer $\mathbf{u}_{0}.$ Moreover, if we set
\begin{align}
e_{0}  & =\ E(\mathbf{u}_{0});\ \ p_{0}=\ C(\mathbf{u}_{0})\label{2}\\
\Gamma & =\Gamma\left(  e_{0},p_{0}\right)  =\left\{  \mathbf{u}\in
X\ |\ E(\mathbf{u})=e_{0},\ C(\mathbf{u})=p_{0}\right\}  ,\label{3}%
\end{align}
every $\mathbf{u}\in\Gamma\ $ is a hylomorphic soliton according to definition
\ref{tdc}.
\end{theorem}

\textbf{Proof}. The proof of this theorem is in \cite{befo2011}. Here we just
give an idea of it. Let $\mathbf{u}_{n}$ be a minimizing sequence of $J. $ $J$
is $G$-compact, then, for a suitable subsequence $\mathbf{u}_{n_{k}} $ and a
suitable sequence $g_{k}$, we get $g_{k}\mathbf{u}_{n_{k}}\rightarrow
\mathbf{u}_{0}.$ Clearly $\mathbf{u}_{0}$ is a minimizer of $J.$

Now let $\Gamma$ be defined as in (\ref{3}). It remains to show that every
$\mathbf{u}\in\Gamma\ $ is a hylomorphic soliton according to definition
\ref{tdc}. First of all notice that $\mathbf{u}_{0}$ is a minimizer of $E$ on
the set
\[
\mathfrak{M}_{p_{0}}=\left\{  \mathbf{u}\in X\ |\ C(\mathbf{u})=p_{0}\right\}
\]
and hence, according to definition \ref{tdc}, every $\mathbf{u}\in\Gamma$ is a
hylomorphic soliton provided that $\Gamma$ satisfies (i), (ii), (iii) of
definition \ref{ds}. Clearly (i) and (iii) are satisfied. In order to prove
(ii), namely that $\Gamma$ is stable, we set
\begin{equation}
V\left(  \mathbf{u}\right)  =\left(  E\left(  \mathbf{u}\right)
-e_{0}\right)  ^{2}+\left(  C\left(  \mathbf{u}\right)  -c_{0}\right)
^{2}.\label{liap}%
\end{equation}
It can be shown that $V$ is a Liapunov function. Then it is sufficient to
apply the classical Liapunov theorem.

$\square$

\bigskip

\begin{remark}
The reader may wonder why we use the functional $J$ rather than mimimizing $E
$ on the manifold $\mathfrak{M}_{p}$, $p\in\mathbb{R}$. As matter of fact, in
general $E$ does not have a minimum on $\mathfrak{M}_{p};$ on the contrary, if
you choose $p_{0}$ given by (\ref{2}), $E$ has a minimum on $\mathfrak{M}%
_{p_{0}}.$ In general, there is a set $I$ of real values such that $\delta\in
I$ implies that $J$ given by (\ref{jd}) is $G$-compact; then for every
$\delta\in I$, there is a $p=p(\delta)$ such that $E$ has a minimum on
$\mathfrak{M}_{p(\delta)}.$Moreover, if you perform a numerical simulations,
it is more efficient to minimize the functional $J$ rather than the functional
$E$ constrained on $\mathfrak{M}_{p(\delta)}$
\end{remark}

\bigskip

\section{The existence result}

\subsection{Statement of the main results\label{smr}}

Equation (\ref{1}) has a variational structure, namely it is the
Euler-Lagrange equation with respect to the functional%
\begin{equation}
S=\frac{1}{2}\int\int\left(  u_{t}^{2}-u_{xx}^{2}\right)  dx\ dt-\int\int
W(u)dx\ dt.\label{action}%
\end{equation}
The Lagrangian relative to the action (\ref{action}) is%
\begin{equation}
\mathcal{L=}\frac{1}{2}\left(  u_{t}^{2}-u_{xx}^{2}\right)
-W(u).\label{lagrangian}%
\end{equation}
This Lagrangian does not depend on $t$ and $x.$ Then, by Noether's Theorem
(see e.g. \cite{gel}, \cite{befogranas}), the energy $E$ and the momentum $C$
defined by%
\[
E=\int\left(  \frac{\partial\mathcal{L}}{\partial u_{t}}u_{t}-\mathcal{L}%
\right)  dx=\frac{1}{2}\int\left(  u_{t}^{2}+u_{xx}^{2}\right)  dx+\int W(u)dx
\]%
\[
C=-\int\left(  \frac{\partial\mathcal{L}}{\partial u_{t}}u_{x}\right)
dx=-\int u_{t}u_{x}\ dx
\]
are constant along the solutions of (\ref{1}).

Equation (\ref{1}), can be rewritten as an Hamiltonian system as follows:%
\begin{equation}
\left\{
\begin{array}
[c]{c}%
\partial_{t}u=v\\
\\
\partial_{t}v=-\partial_{x}^{4}u-W^{\prime}(u)
\end{array}
\right. \label{primo}%
\end{equation}

The phase space is given by%
\[
X=H^{2}(\mathbb{R})\times L^{2}(\mathbb{R})
\]
and the generic point in $X$ will be denoted by%
\[
\mathbf{u}=\left[
\begin{array}
[c]{c}%
u\\
v
\end{array}
\right]  .
\]
Here $H^{2}(\mathbb{R})$ denotes the usual Sobolev space.

The norm of $X$ is given by%
\[
\left\Vert \mathbf{u}\right\Vert =\left(  \int\left(  v^{2}+u_{xx}^{2}%
+u^{2}\right)  dx\right)  ^{\frac{1}{2}}.
\]
The energy and the momentum, as functionals defined on $X,$ take the following
form%
\[
E\left(  \mathbf{u}\right)  =\frac{1}{2}\int\left(  v^{2}+u_{xx}^{2}\right)
dx+\int W(u)dx
\]%
\[
C\left(  \mathbf{u}\right)  =-\int vu_{x}\ dx.
\]
\bigskip

Next, we will apply the abstract theory of section \ref{HS} where the momentum
$C\left(  \mathbf{u}\right)  $ plays the role of the hylenic charge.

We make the following assumptions:

\begin{itemize}
\item (W-i) \textbf{(Positivity}) $\exists\eta>0\ $such that $W(s)\geq\eta
s^{2}$ for $|s|\leq1$ and $W(s)\geq\eta$ for $|s|\geq1.$

\item (W-ii) \textbf{(Nondegeneracy at 0})$\ W^{\prime\prime}(0)=1$

\item (W-iii) \textbf{(Hylomorphy})\ $\exists M>0,\ \exists\alpha\in
\lbrack0,2),\forall s\geq0,$%
\[
W(s)\leq M\left\vert s\right\vert ^{\alpha}.
\]

\end{itemize}

Here there are some comments on assumptions (W-ii),(W-iii).

(W-ii) The assumption $W^{\prime\prime}(0)=1$ can be weakened just assuming
the existence of$\ W^{\prime\prime}(0).$In fact, by (W-i) we have
$\;W^{\prime\prime}(0)>0$ and we can reduce to the case $W^{\prime\prime
}(0)=1,$ by rescaling space and time. By this assumption we can write%
\begin{equation}
W(s)=\frac{1}{2}s^{2}+N(s),\ \ N(s)=o(s^{2}).\label{NN}%
\end{equation}

(W-iii) This is the crucial assumption which characterizes the potentials
which might produce hylomorphic solitons; notice that this assumptions
concerns $W$ only for the positive values of $s.$

\bigskip

We have the following results:

\begin{theorem}
\label{main-theorem} Assume that (W-i),(W-ii),(W-iii) hold, then there exists
an open interval $I\subset\mathbb{R}$ such that, for every $\delta\in I,$
there is an hylomorphic soliton $\mathbf{u}_{\delta}$ for the dynamical system
(\ref{primo}) . Moreover, if $\delta_{1}\neq\delta_{2},$ $\mathbf{u}%
_{\delta_{1}}\neq g\mathbf{u}_{\delta_{2}}$ for every $g\in G$.
\end{theorem}

\begin{theorem}
\label{imp}Let $\mathbf{u}_{\delta}=\left(  u_{\delta},v_{\delta}\right)  $ be
a soliton as in Theorem \ref{main-theorem}. Then the solution of
eq.(\ref{1})$\ $with initial data $\left(  u_{\delta},v_{\delta}\right)  $ has
the following form:
\[
u(t,x)=u_{\delta}(x-ct)
\]
where $u_{\delta}\ $is a solution of the following equation
\begin{equation}
\frac{\partial^{4}u_{\delta}}{\partial x^{4}}+c^{2}\frac{\partial^{2}%
u_{\delta}}{\partial x^{2}}+W^{\prime}(u_{\delta})=0\label{stat2}%
\end{equation}
and $c$ is a constant which depends on $u_{\delta}.$
\end{theorem}

\bigskip

\begin{remark}
So we get the existence of solutions of (\ref{stat2}) by a different proof
from that in \cite{McKW90} and \cite{SaWei09}. We point out that (\ref{stat2})
could have solutions which are not minimizers. In this case these solutions
give rise to solitary waves which are not solitons.
\end{remark}

The proofs of Theorem \ref{main-theorem} and of Theorem \ref{imp} will be
given in the next section.

\subsection{Proof of the main results}

By (W-ii), we have that for $\mathbf{u}=\left[
\begin{array}
[c]{c}%
u\\
v
\end{array}
\right]  \in X=H^{2}(\mathbb{R})\times L^{2}(\mathbb{R})$%
\begin{equation}
E\left(  \mathbf{u}\right)  =\frac{1}{2}\left\Vert \mathbf{u}\right\Vert
^{2}+\int N(u)dx.\label{denot}%
\end{equation}

\begin{lemma}
\label{zero}Let $M>0.$ Then there exists a constant $C>0$ such that $\left(
E\left(  \mathbf{u}\right)  \leq M\right)  \Rightarrow$ $\left(  \left\Vert
\mathbf{u}\right\Vert \leq C\right)  $.
\end{lemma}

\textbf{Proof. }Assume that
\begin{equation}
E\left(  \mathbf{u}\right)  =\frac{1}{2}\int\left(  v^{2}+u_{xx}^{2}\right)
dx+\int W(u)dx\leq M.\label{laida}%
\end{equation}
Then, since $\ W(u)\geq0,$ we have that
\begin{equation}
\int\left(  v^{2}+u_{xx}^{2}\right)  dx\leq M.\label{see}%
\end{equation}

It remains to prove that also%
\begin{equation}
\int u^{2}dx\text{ is bounded.}\label{bou}%
\end{equation}

We now set
\[
\Omega_{u}^{+}=\left\{  x\ |\ u(x)>1\right\}  ;\ \Omega_{u}^{-}=\left\{
x\ |\ u(x)<-1\right\}  .
\]
Then, if (\ref{laida}) holds, by (W-i) we have%
\begin{equation}
M\geq\int W(u)dx^{+}\geq\int_{\Omega_{u}^{+}\cup\Omega_{u}^{-}}W(u)dx\geq
\eta\left\vert \Omega_{u}^{+}\right\vert +\eta\left\vert \Omega_{u}%
^{-}\right\vert \label{ba}%
\end{equation}
where $\left\vert \Omega\right\vert $ denotes the measure of $\Omega.$ Now we
show that
\begin{equation}
\int_{\Omega_{u}^{+}}u^{2}dx\text{ is bounded.}\label{lalla}%
\end{equation}
Set $v=u-1,$ then, since $v=0$ on $\partial\Omega_{u}^{+},$ by the
Poincar\`{e} inequality, there is a constant $c>0$ such that
\begin{equation}
\int_{\Omega_{u}^{+}}v^{2}dx\leq c\int_{\Omega_{u}^{+}}v_{x}^{2}%
dx.\label{lulla}%
\end{equation}
since we are in dimension one, it is easy to check that $c\leq\left\vert
\Omega_{u}^{+}\right\vert ^{2}.$

On the other hand
\begin{equation}
\int_{\Omega_{u}^{+}}v_{x}^{2}dx=-\int_{\Omega_{u}^{+}}v\text{ }v_{xx}%
dx\leq\left\Vert v\right\Vert _{L^{2}(\Omega_{u}^{+})}\left\Vert
v_{xx}\right\Vert _{L^{2}(\Omega_{u}^{+})}.\label{pa}%
\end{equation}

Then, since $v=u-1,$ by (\ref{lulla}) and (\ref{pa}),
\[
\left\Vert u-1\right\Vert _{L^{2}(\Omega_{u}^{+})}^{2}\leq c\left\Vert
u-1\right\Vert _{L^{2}(\Omega_{u}^{+})}\left\Vert u_{xx}\right\Vert
_{L^{2}(\Omega_{u}^{+})}%
\]

we easily get%
\begin{equation}
\left\Vert u\right\Vert _{L^{2}(\Omega_{u}^{+})}^{2}-2\left\vert \Omega
_{u}^{+}\right\vert ^{\frac{1}{2}}\left\Vert u\right\Vert _{L^{2}(\Omega
_{u}^{+})}+\left\vert \Omega_{u}^{+}\right\vert \leq c\left(  \left\Vert
u\right\Vert _{L^{2}(\Omega_{u}^{+})}+\left\vert \Omega_{u}^{+}\right\vert
\right)  \left\Vert u_{xx}\right\Vert _{L^{2}(\Omega_{u}^{+})}.\label{sue}%
\end{equation}

By (\ref{see}) and (\ref{ba}) we have
\begin{equation}
\left\Vert u_{xx}\right\Vert _{L^{2}(\Omega_{u}^{+})}\leq\sqrt{M},\text{
}\left\vert \Omega_{u}^{+}\right\vert \leq\frac{M}{\eta}.\label{pe}%
\end{equation}

By (\ref{sue}) and (\ref{pe}) we get
\[
\left\Vert u\right\Vert _{L^{2}(\Omega_{u}^{+})}^{2}-2\left(  \frac{M}{\eta
}\right)  ^{\frac{1}{2}}\left\Vert u\right\Vert _{L^{2}(\Omega_{u}^{+})}\leq
c\sqrt{M}\left(  \left\Vert u\right\Vert _{L^{2}(\Omega_{u}^{+})}+\frac
{M}{\eta}\right)  .
\]
From which we easily deduce (\ref{lalla}). Analogously, w get also that
\begin{equation}
\int_{\Omega_{u}^{-}}u^{2}dx\text{ is bounded.}\label{lalla-}%
\end{equation}

By (W-i)
\[
M\geq\int W(u)dx=\int_{\left\vert u(x)\right\vert \leq1}W(u(x))dx+\int
_{\Omega_{u}^{+}\cup\Omega_{u}^{-}}W(u(x))dx\geq\eta\int_{|u(x)|\leq1}u^{2}dx.
\]
So, by (\ref{lalla}), (\ref{lalla-}) and the above inequality, there is a
constant $R$ such that%
\[
\int u^{2}dx=\int_{|u(x)|\leq1}u^{2}dx+\int_{\Omega_{u}^{+}\cup\Omega_{u}^{-}%
}u^{2}dx\leq\frac{M}{\eta}+R.
\]
We conclude that $\int u^{2}dx$ is bounded.

$\square$

\bigskip

\begin{lemma}
\label{zeropiu}Let $\mathbf{u}_{n}$ be a sequence in $X$ such that%
\begin{equation}
E\left(  \mathbf{u}_{n}\right)  \rightarrow0.\label{atto}%
\end{equation}
Then, up to a subsequence, we have $\left\Vert \mathbf{u}_{n}\right\Vert
_{X}\rightarrow0.$
\end{lemma}

\textbf{Proof.} Let $\mathbf{u}_{n}=(u_{n},v_{n}),$ $u_{n}\in H^{2}%
(\mathbb{R)},v_{n}\in L^{2}(\mathbb{R)},$ be a sequence such that $E\left(
\mathbf{u}_{n}\right)  \rightarrow0.$ Then clearly $\left\Vert v_{n}%
\right\Vert _{L^{2}}\rightarrow0.$ By Lemma \ref{zero}, $u_{n}$ is bounded in
$H^{2}(\mathbb{R})$ and hence, by the Sobolev embedding therems, $u_{n}$ is
bounded in $L^{\infty}(\mathbb{R)}$, moreover for all $n$ we have
$u_{n}(x)\rightarrow0$ for $\left\vert x\right\vert \rightarrow\infty.$

For each $n$ let $\tau_{n}$ be a maximum point of $\left\vert u_{n}\right\vert
$ and set%
\[
u_{n}^{\prime}(x)=u_{n}(\tau_{n}+x),\text{ }v_{n}^{\prime}(x)=v_{n}(\tau
_{n}+x),
\]
so that
\begin{equation}
\left\vert u_{n}^{\prime}(0)\right\vert =\max\left\vert u_{n}^{\prime
}\right\vert .\label{prov}%
\end{equation}

$.$

Clearly $u_{n}^{\prime}$ is bounded in $H^{2}(\mathbb{R}),$ then, up to a
subsequence, we get%
\begin{equation}
u_{n}^{\prime}\rightharpoonup u\text{ weakly in }H^{2}(\mathbb{R})\label{a}%
\end{equation}
and consequently%
\begin{equation}
\frac{d^{2}u_{n}^{\prime}}{dx^{2}}\rightharpoonup\frac{d^{2}u}{dx^{2}}\text{
weakly in }L^{2}(\mathbb{R)}\text{.}\label{b}%
\end{equation}
On the other end, since $E\left(  \mathbf{u}_{n}\right)  \rightarrow0,$ we
have $\frac{d^{2}u_{n}}{dx^{2}}\rightarrow0$ in $L^{2}(\mathbb{R)}$. Then
also
\begin{equation}
\frac{d^{2}u_{n}^{\prime}}{dx^{2}}\rightarrow0\text{ in }L^{2}(\mathbb{R)}%
\text{.}\label{c}%
\end{equation}
From (\ref{b}) and (\ref{c}) we get
\[
\frac{d^{2}u}{dx^{2}}=0.
\]
So $u\in$ $H^{2}(\mathbb{R})$ is linear and consequently%
\begin{equation}
u=0.\label{d}%
\end{equation}

Now set
\[
B_{R}=\left\{  x\in\mathbb{R}:\left\vert x\right\vert <R\right\}  ,\text{ }R>0
\]
then, by the compact embedding $H^{2}(B_{R})\subset\subset L^{\infty}(B_{R}),$
by (\ref{a}) and (\ref{d}), we get
\begin{equation}
u_{n}^{\prime}\rightarrow0\text{ in }L^{\infty}(B_{R}).\label{e}%
\end{equation}
By (\ref{prov}) and (\ref{e}) we get%
\[
\left\Vert u_{n}^{\prime}\right\Vert _{L^{\infty}(\mathbb{R})}=\left\vert
u_{n}^{\prime}(0)\right\vert \rightarrow0.
\]
So, if $n$ is sufficiently large, we have $\left\vert u_{n}^{\prime
}(x)\right\vert \leq1$ for all $x.$

Then, setting $\mathbf{u}_{n}^{\prime}=(u_{n}^{\prime},v_{n}^{\prime}),$ by
(W-i), we have that
\begin{align}
E\left(  \mathbf{u}_{n}^{\prime}\right)   & =\int(\frac{1}{2}\left(
v_{n}^{\prime2}+\left(  \partial_{xx}^{2}u_{n}^{\prime}\right)  ^{2}%
)+W(u_{n}^{\prime})\right)  dx\nonumber\\
& \geq\int\left(  \frac{1}{2}\left(  v_{n}^{\prime2}+\left(  \partial_{xx}%
^{2}u_{n}^{\prime}\right)  ^{2})\right)  +\eta u_{n}^{\prime2}\right)
dx\nonumber\\
& \geq c\left\Vert \mathbf{u}_{n}^{\prime}\right\Vert ^{2}\label{att}%
\end{align}
where $c$ is a positive constant.

Since%
\[
E\left(  \mathbf{u}_{n}^{\prime}\right)  =E\left(  \mathbf{u}_{n}\right)
,\left\Vert \mathbf{u}_{n}^{\prime}\right\Vert =\left\Vert \mathbf{u}%
_{n}\right\Vert ,
\]
by (\ref{att}), (\ref{atto}) we have
\[
\left\Vert \mathbf{u}_{n}\right\Vert _{X}\rightarrow0\text{ .}%
\]

$\square$

We set
\begin{align*}
\Lambda_{0}  & =\ \underset{\mathbf{u}\in X}{\inf}\ \frac{\frac{1}%
{2}\left\Vert \mathbf{u}\right\Vert ^{2}}{\left\vert C\left(  \mathbf{u}%
\right)  \right\vert },\text{ }\\
\Lambda_{\ast}  & =\ \underset{\mathbf{u}\in X}{\inf}\ \frac{E\left(
\mathbf{u}\right)  }{\left\vert C\left(  \mathbf{u}\right)  \right\vert
}=\ \underset{\mathbf{u}\in X}{\inf}\ \frac{\frac{1}{2}\left\Vert
\mathbf{u}\right\Vert ^{2}+\int N(u)dx}{\left\vert C\left(  \mathbf{u}\right)
\right\vert }.
\end{align*}

\begin{lemma}
\label{uno}The following inequality holds:%
\[
\Lambda_{0}\geq1.
\]

\end{lemma}

\textbf{Proof}: For $\mathbf{u}=(v,u)$ we have%
\begin{align*}
\left\vert C\left(  \mathbf{u}\right)  \right\vert  & \leq\int\left\vert
v\partial_{x}u\ \right\vert \ dx\leq\left(  \int v^{2}\ dx\right)  ^{1/2}%
\cdot\left(  \int\left\vert \partial_{x}u\right\vert ^{2}\ dx\right)  ^{1/2}\\
& \leq\frac{1}{2}\int v^{2}\ dx+\frac{1}{2}\int\left\vert \partial
_{x}u\right\vert ^{2}\ dx\\
& =\frac{1}{2}\int v^{2}\ dx-\frac{1}{2}\int uu_{xx}\ dx\\
& \leq\frac{1}{2}\int v^{2}\ dx+\frac{1}{2}\int\frac{1}{2}\left[  u^{2}%
+u_{xx}^{2}\right]  \ dx\\
& \leq\frac{1}{2}\int\left[  v^{2}\ +u_{xx}^{2}+u^{2}\right]  \ dx=\frac{1}%
{2}\left\Vert \mathbf{u}\right\Vert ^{2}.
\end{align*}
Then, for every $\mathbf{u}$
\[
\Lambda_{0}\geq\ \frac{\frac{1}{2}\left\Vert \mathbf{u}\right\Vert ^{2}%
}{\left\vert C\left(  \mathbf{u}\right)  \right\vert }\geq1.
\]

$\square$

\bigskip

The next lemma provides a crucial estimate for the existence of solitons:

\begin{lemma}
\label{due}We have%
\[
\Lambda_{\ast}<1
\]

\end{lemma}

\textbf{Proof:}\ Let $U\in C^{2}$ be a positive function with compact support
such that%
\begin{equation}
\frac{\int\left(  U_{xx}\right)  ^{2}}{\int\left(  U_{x}\right)  ^{2}}%
<\frac{1}{2}.\label{UU}%
\end{equation}
Such a function exists; in fact if $U_{0}$ is any positive function with
compact support, $U(x)=U_{0}\left(  \frac{x}{\lambda}\right)  $ satisfies
(\ref{UU}) for $\lambda$ sufficiently large. Take
\[
\mathbf{u}_{R}=\left(  u_{R},v\right)  =\left(  RU,RU_{x}\right)  .
\]
By the definition of $X$, $\mathbf{u}_{R}\in X.$ Now we can estimate
$\Lambda_{\ast}$:%

\begin{align*}
\Lambda_{\ast} &  =\ \underset{\mathbf{u}\in X}{\inf}\ \frac{\frac{1}%
{2}\left\Vert \mathbf{u}\right\Vert ^{2}+\int N(u)dx}{\left\vert C\left(
\mathbf{u}\right)  \right\vert }\leq\frac{\frac{1}{2}\left\Vert \mathbf{u}%
_{R}\right\Vert ^{2}+\int N(u_{R})dx}{\left\vert C\left(  \mathbf{u}%
_{R}\right)  \right\vert }\\
&  =\frac{\frac{1}{2}\int\left[  \left(  RU_{x}\right)  ^{2}+\left(
RU_{xx}\right)  ^{2}+\left(  RU\right)  ^{2}\right]  dx+\int N(RU)dx}%
{\int\left(  RU_{x}\right)  ^{2}\ dx}\\
&  =\frac{\frac{1}{2}\int\left[  \left(  RU_{x}\right)  ^{2}+\left(
RU_{xx}\right)  ^{2}\right]  dx}{\int\left(  RU_{x}\right)  ^{2}\ dx}%
+\frac{\int W(RU)dx}{\int\left(  RU_{x}\right)  ^{2}\ dx}\\
&  =\frac{1}{2}+\frac{1}{2}\frac{\int\left(  U_{xx}\right)  ^{2}dx}%
{\int\left(  U_{x}\right)  ^{2}\ dx}+\frac{\int W(RU)dx}{\int\left(
RU_{x}\right)  ^{2}\ dx}\text{ }\ \text{(by (W-iii))}\\
&  \leq\frac{1}{2}+\frac{1}{2}\frac{\int\left(  U_{xx}\right)  ^{2}dx}%
{\int\left(  U_{x}\right)  ^{2}\ dx}+\frac{\int M\left\vert RU\right\vert
^{\alpha}dx}{\int\left(  RU_{x}\right)  ^{2}\ dx}\ \ \text{(by (\ref{UU}))}\\
&  <\frac{1}{2}+\frac{1}{4}+\frac{M}{R^{2-\alpha}}\cdot\frac{\int\left\vert
U\right\vert ^{\alpha}dx}{\int U_{x}^{2}\ dx}.
\end{align*}
Then, for $R$ sufficiently large, we get the conclusion.

$\square$

\begin{lemma}
\label{tre}Consider any sequence
\[
\mathbf{u}_{n}=\mathbf{u}+\mathbf{w}_{n}\in X
\]
where $\mathbf{w}_{n}$ converges weakly to $0.$ Then%
\begin{equation}
E(\mathbf{u}_{n})=E(\mathbf{u})+E(\mathbf{w}_{n})+o(1)\label{mi}%
\end{equation}
and%
\begin{equation}
C(\mathbf{u}_{n})=C(\mathbf{u})+C(\mathbf{w}_{n})+o(1).\label{mu}%
\end{equation}

\end{lemma}

\textbf{Proof. }First of all we introduce the following notation:%
\[
K(u)=\int N\left(  u\right)  dx\text{ and }K_{\Omega}(u)=\int_{\Omega}N\left(
u\right)  dx\text{, }\Omega\text{ open subset in }\mathbb{R}\text{.}%
\]

As usual $u,w_{n}$ will denote the first components respectively of
$\mathbf{u,w}_{n}\in H^{2}\left(  \mathbb{R}\right)  \times L^{2}\left(
\mathbb{R}\right)  $.

We have to show that $\underset{n\rightarrow\infty}{\lim}\left\vert E\left(
\mathbf{u}+\mathbf{w}_{n}\right)  -E\left(  \mathbf{u}\right)  -E\left(
\mathbf{w}_{n}\right)  \right\vert =0.$ By ( \ref{denot}) we have that%
\begin{align}
&  \underset{n\rightarrow\infty}{\lim}\left\vert E\left(  \mathbf{u}%
+\mathbf{w}_{n}\right)  -E\left(  \mathbf{u}\right)  -E\left(  \mathbf{w}%
_{n}\right)  \right\vert \label{inf}\\
&  \leq\ \underset{n\rightarrow\infty}{\lim\frac{1}{2}}\left\vert \left\Vert
\mathbf{u}+\mathbf{w}_{n}\right\Vert ^{2}-\left\Vert \mathbf{u}\right\Vert
^{2}-\left\Vert \mathbf{w}_{n}\right\Vert ^{2}\right\vert \nonumber\\
&  +\ \underset{n\rightarrow\infty}{\lim}\left\vert \int\left(  N\left(
u+w_{n}\right)  -N\left(  u\right)  -N\left(  w_{n}\right)  \right)
dx\right\vert .\nonumber
\end{align}

If $(\cdot,\cdot)$ denotes the inner product induced by the norm $\left\Vert
\cdot\right\Vert $ we have:%
\begin{equation}
\underset{n\rightarrow\infty}{\lim}\left\vert \left\Vert \mathbf{u}%
+\mathbf{w}_{n}\right\Vert ^{2}-\left\Vert \mathbf{u}\right\Vert
^{2}-\left\Vert \mathbf{w}_{n}\right\Vert ^{2}\right\vert =\ \underset
{n\rightarrow\infty}{\lim}\left\vert 2\left(  \mathbf{u},\mathbf{w}%
_{n}\right)  \right\vert =0.\label{min}%
\end{equation}

Then by (\ref{inf}) and (\ref{min}) we have%
\begin{align}
& \underset{n\rightarrow\infty}{\lim}\left\vert E\left(  \mathbf{u}%
+\mathbf{w}_{n}\right)  -E\left(  \mathbf{u}\right)  -E\left(  \mathbf{w}%
_{n}\right)  \right\vert \label{sotto}\\
& \leq\ \underset{n\rightarrow\infty}{\lim}\left\vert \int\left(  N\left(
u+w_{n}\right)  -N\left(  u\right)  -N\left(  w_{n}\right)  \right)
dx\right\vert .
\end{align}

Choose $\varepsilon>0$ and $R=R(\varepsilon)>0$ such that
\begin{equation}
\left\vert \int_{B_{R}^{c}}N\left(  u\right)  \right\vert <\varepsilon,\text{
}\int_{B_{R}^{c}}\left\vert u\right\vert <\varepsilon\label{bis}%
\end{equation}

$\ $where $\ $%
\[
B_{R}^{c}=\mathbb{R}^{N}-B_{R}\text{ and }B_{R}=\left\{  x\in\mathbb{R}%
^{N}:\left\vert x\right\vert <R\right\}  .
\]
Since $w_{n}\rightharpoonup0$ weakly in $H^{2}\left(  \mathbb{R}\right)  $, by
usual compactness arguments, we have that
\begin{equation}
K_{B_{R}}\left(  w_{n}\right)  \rightarrow0\text{ and }K_{B_{R}}\left(
u+w_{n}\right)  \rightarrow K_{B_{R}}\left(  u\right)  .\label{bibis}%
\end{equation}
Then, by (\ref{bis}) and (\ref{bibis}), we have%

\begin{align}
&  \underset{n\rightarrow\infty}{\lim}\left\vert \int\left[  N\left(
u+w_{n}\right)  -N\left(  u\right)  -N\left(  w_{n}\right)  \right]
\right\vert \nonumber\\
&  =\ \underset{n\rightarrow\infty}{\lim}|K_{B_{R}^{c}}\left(  u+w_{n}\right)
+K_{B_{R}}\left(  u+w_{n}\right) \nonumber\\
&  \ \ \ \ \ \ \ \ \ \ \ \ \ \ -K_{B_{R}^{c}}\left(  u\right)  -K_{B_{R}%
}\left(  u\right)  -K_{B_{R}^{c}}\left(  w_{n}\right)  -K_{B_{R}}\left(
w_{n}\right)  |
\end{align}
Then, by (\ref{bibis}) and (\ref{bis})
\begin{align*}
&  \underset{n\rightarrow\infty}{\lim}\left\vert \int\left[  N\left(
u+w_{n}\right)  -N\left(  u\right)  -N\left(  w_{n}\right)  \right]
\right\vert \\
&  =\ \underset{n\rightarrow\infty}{\lim}\left\vert K_{B_{R}^{c}}\left(
u+w_{n}\right)  -K_{B_{R}^{c}}\left(  u\right)  -K_{B_{R}^{c}}\left(
w_{n}\right)  \right\vert \\
&  \leq\ \underset{n\rightarrow\infty}{\lim}\left\vert K_{B_{R}^{c}}\left(
u+w_{n}\right)  -K_{B_{R}^{c}}\left(  w_{n}\right)  \right\vert +\varepsilon.
\end{align*}

By the intermediate value theorem there are $\zeta_{n}\ $in $(0,1)$ such that
\begin{equation}
\left\vert K_{B_{R}^{c}}\left(  u+w_{n}\right)  -K_{B_{R}^{c}}\left(
w_{n}\right)  \right\vert =\int_{B_{R}^{c}}N^{\prime}\left(  \zeta_{n}%
u+w_{n}\right)  udx.\label{fere}%
\end{equation}
Since $w_{n}$ is bounded in $H^{2}\left(  \mathbb{R}\right)  ,$ $\zeta
_{n}u+w_{n}$ is bounded in $L^{\infty},$ so that there exists a positive
constant $M$ such that
\begin{equation}
\left\Vert N^{\prime}\left(  \zeta_{n}u+w_{n}\right)  \right\Vert _{L^{\infty
}}\leq M.\label{feri}%
\end{equation}
By (\ref{fere}), (\ref{feri}) and (\ref{bis}) we have
\begin{equation}
\left\vert K_{B_{R}^{c}}\left(  u+w_{n}\right)  -K_{B_{R}^{c}}\left(
w_{n}\right)  \right\vert \leq M\int_{B_{R}^{c}}\left\vert u\right\vert
<M\varepsilon.\label{get}%
\end{equation}

Then, by (\ref{pepe}) and (\ref{get}), we get%
\begin{equation}
\underset{n\rightarrow\infty}{\lim}\left\vert \int\left[  N\left(
u+w_{n}\right)  -N\left(  u\right)  -N\left(  w_{n}\right)  \right]
\right\vert \leq\varepsilon+M\cdot\varepsilon.\label{sopra}%
\end{equation}
Finally by (\ref{sotto}) and (\ref{sopra}) and since $\varepsilon$ is arbitray
we get
\[
\underset{n\rightarrow\infty}{\lim}\left\vert E\left(  \mathbf{u}%
+\mathbf{w}_{n}\right)  -E\left(  \mathbf{u}\right)  -E\left(  \mathbf{w}%
_{n}\right)  \right\vert =0
\]

and so (\ref{mi}) is proved$.$ The proof of (\ref{mu}) is immediate.

$\square$

\bigskip

By lemma \ref{uno} and lemma \ref{due}, we have that
\[
\Lambda_{\ast}<\Lambda_{0}.
\]
So there exist $\mathbf{u}_{0}\in X$ and $b>0$ such that
\[
\frac{E(\mathbf{u}_{0})}{\left\vert C(\mathbf{u}_{0})\right\vert }\leq
\Lambda_{0}-b.
\]
Then we can choose $\delta>0$ such that%
\begin{equation}
\frac{E(\mathbf{u}_{0})}{\left\vert C(\mathbf{u}_{0})\right\vert }+\delta
E(\mathbf{u}_{0})\leq\Lambda_{0}-\frac{b}{2}\label{kkk}%
\end{equation}
and we define
\begin{equation}
J(\mathbf{u})=\frac{E(\mathbf{u})}{\left\vert C(\mathbf{u})\right\vert
}+\delta E(\mathbf{u}).\label{J}%
\end{equation}

Then we have that%
\begin{equation}
J_{\ast}:=\inf J(\mathbf{u})\leq J(\mathbf{u}_{0})\leq\Lambda_{0}-\frac{b}%
{2}.\label{brutta}%
\end{equation}
\bigskip

\begin{lemma}
\label{GC}The functional defined by (\ref{J}) is $G$-compact (where $G$ is
defined by (\ref{ggg})).
\end{lemma}

\textbf{Proof. }Let $\mathbf{u}_{n}$ $=\left(  u_{n},v_{n}\right)  $ be a
minimizing sequence for $J.$ Since the $G$-compactness depends on
subsequences, we can take a subsequence in which all the $C(\mathbf{u}_{n})$
have the same sign. So, to fix the ideas, we can assume that
\begin{equation}
C(\mathbf{u}_{n})>0;\label{piu}%
\end{equation}
thus we have that
\[
J(\mathbf{u}_{n})=\frac{E(\mathbf{u}_{n})}{C(\mathbf{u}_{n})}+\delta
E(\mathbf{u}_{n}).
\]
It is immediate to see that $E(\mathbf{u}_{n})=\frac{1}{2}\left\Vert
\mathbf{u}_{n}\right\Vert ^{2}+\int N(u_{n})dx$ is bounded. Then, by lemma
\ref{zero}, $\left\Vert \mathbf{u}_{n}\right\Vert $ is bounded and hence,
passing eventually to a suitable subsequence, we have $\mathbf{u}%
_{n}\rightharpoonup\mathbf{u}$ weakly in $X.$ Now, starting from
$\mathbf{u}_{n},$ we construct a minimizing sequence $\mathbf{u}_{n}^{\prime}$
which weakly converges to
\begin{equation}
\mathbf{\bar{u}\neq}0.\label{ss}%
\end{equation}

To this end we first show that:%

\begin{equation}
\left\Vert u_{n}\right\Vert _{L^{\infty}}\text{ does not converge to
}0.\label{ab}%
\end{equation}
Arguing by contradiction, assume that
\[
\left\Vert u_{n}\right\Vert _{L^{\infty}}\text{ }\rightarrow0.
\]
Then, since $N(s)=o(s^{2}),$ there is a sequence of positive real numbers
$\varepsilon_{n}$ with $\varepsilon_{n}\rightarrow0$ such that
\[
\frac{E(\mathbf{u}_{n})}{C(\mathbf{u}_{n})}\geq\frac{\frac{1}{2}\left(
\left\Vert \frac{d^{2}u_{n}}{dx^{2}}\right\Vert _{L^{2}}^{2}+\left\Vert
u_{n}\right\Vert _{L^{2}}^{2}\right)  -\int\left\vert N(u_{n})\right\vert
dx}{C(\mathbf{u}_{n})}%
\]%
\[
=\frac{\frac{1}{2}\left(  \left\Vert \frac{d^{2}u_{n}}{dx^{2}}\right\Vert
_{L^{2}}^{2}+\left\Vert u_{n}\right\Vert _{L^{2}}^{2}\right)  -\frac
{\varepsilon_{n}}{2}\left\Vert u_{n}\right\Vert _{L^{2}}^{2}}{C(\mathbf{u}%
_{n})}\geq
\]%
\[
\geq\frac{\frac{1}{2}\left(  \left\Vert \frac{d^{2}u_{n}}{dx^{2}}\right\Vert
_{L^{2}}^{2}+\left\Vert u_{n}\right\Vert _{L^{2}}^{2}\right)  -\frac
{\varepsilon_{n}}{2}\left(  \left\Vert \frac{d^{2}u_{n}}{dx^{2}}\right\Vert
_{L^{2}}^{2}+\left\Vert u_{n}\right\Vert _{L^{2}}^{2}\right)  }{C(\mathbf{u}%
_{n})}=
\]%
\[
\frac{\frac{1}{2}\left(  \left\Vert \frac{d^{2}u_{n}}{dx^{2}}\right\Vert
_{L^{2}}^{2}+\left\Vert u_{n}\right\Vert _{L^{2}}^{2}\right)  }{C(\mathbf{u}%
_{n})}(1-\varepsilon_{n})\geq(\text{by definition of }\Lambda_{0}\text{)}%
\]%
\[
\geq\Lambda_{0}(1-\varepsilon_{n}).
\]
And hence%
\begin{equation}
J(\mathbf{u}_{n})\geq\Lambda_{0}(1-\varepsilon_{n}).\label{tap}%
\end{equation}

On the other hand by (\ref{brutta})
\begin{equation}
\lim J(\mathbf{u}_{n})\leq\Lambda_{0}-\frac{b}{2}.\label{tip}%
\end{equation}
Clearly (\ref{tip}) contradicts (\ref{tap}).

So (\ref{ab}) holds and consequently there exist $b>0$ and a sequence $x_{n}$
such that, up to a subsequence,
\begin{equation}
\left\vert u_{n}(x_{n})\right\vert \geq b.\label{pippa}%
\end{equation}

Now we set
\[
\mathbf{u}_{n}^{\prime}(x)=\mathbf{u}_{n}(x+x_{n}),\text{ }u_{n}^{\prime
}(x)=u_{n}(x+x_{n}).
\]
Clearly also $\mathbf{u}_{n}^{\prime}(x)$ is a minimizing sequence, moreover,
by (\ref{pippa}),
\begin{equation}
\left\vert u_{n}^{\prime}(0)\right\vert \geq b.\label{poppa}%
\end{equation}
Since, up to a subsequence, $\mathbf{u}_{n}^{\prime}\rightharpoonup
\mathbf{\bar{u}}\in X$ weakly in $X,$ we have, by standard compact embeddings
results, that
\[
u_{n}^{\prime}\rightarrow\bar{u}\text{ in }L^{\infty}(-1,1)\text{ }%
\]
where $\bar{u}$ denotes the first component of $\mathbf{\bar{u}.}$Then by
(\ref{poppa}) we have $\bar{u}\neq0$ and then $\mathbf{\bar{u}\neq}0$. So
(\ref{ss}) is proved.

Now set%

\[
\mathbf{u}_{n}^{\prime}=\mathbf{\bar{u}+\mathbf{w}_{n}}%
\]
with $\mathbf{\mathbf{w}_{n}}\rightharpoonup0$ weakly in $X$.

We finally show that there is no splitting, namely that $\mathbf{\mathbf{w}%
_{n}\rightarrow0}$ strongly in $X.$To this hand first we show that%
\begin{equation}
C(\mathbf{\bar{u}}+\mathbf{\mathbf{w}_{n})}\text{ does not converge to
}0.\label{not}%
\end{equation}
Arguing by contradiction assume that $C(\mathbf{\bar{u}+\mathbf{w}_{n})}$
converges to $0.$ Then, since $\mathbf{\bar{u}+\mathbf{w}_{n}}$ is a
minimizing sequence for $J,$ also $E(\mathbf{\bar{u}+\mathbf{w}_{n})}$
converges to $0$ and then, by Lemma \ref{zeropiu}, we get
\begin{equation}
\mathbf{\bar{u}+\mathbf{w}_{n}}\rightarrow0\text{ in }X.\label{s}%
\end{equation}
From (\ref{s}) and since $\mathbf{\mathbf{w}_{n}}\rightharpoonup0$ weakly in
$X,$ we have that $\mathbf{\bar{u}}=0$, contradicting (\ref{ss}). So
(\ref{not}) holds and, passing enventually to a subsequence, we can assume
\begin{equation}
C(\mathbf{\bar{u}+\mathbf{w}_{n})}\geq\delta>0.\label{far}%
\end{equation}

By lemma \ref{tre}, we have%
\[
E(\mathbf{u}_{n}^{\prime})=E(\mathbf{\bar{u}+\mathbf{w}_{n}})=E(\mathbf{\bar
{u}})+E(\mathbf{w}_{n})+o(1)
\]
and
\begin{equation}
C(\mathbf{u}_{n}^{\prime})=C(\mathbf{\bar{u}+\mathbf{w}_{n})}=C(\mathbf{\bar
{u}})+C(\mathbf{w}_{n})+o(1)\geq(\text{by (\ref{far}))}\geq\delta
>0.\label{fir}%
\end{equation}

Then%
\begin{align*}
J_{\ast} &  :=\lim J(\mathbf{u}_{n}^{\prime})=\lim\frac{E(\mathbf{u}%
_{n}^{\prime})}{C(\mathbf{u}_{n}^{\prime})}+\delta E(\mathbf{u}_{n}^{\prime
})\\
&  =\lim\left[  \frac{E(\mathbf{\bar{u}})+E(\mathbf{w}_{n})+o(1)}%
{C(\mathbf{\bar{u}})+C(\mathbf{w}_{n})+o(1)}+\delta E(\mathbf{\bar{u}})+\delta
E(\mathbf{w}_{n})+o(1)\right] \\
&  =\lim\left[  \frac{E(\mathbf{\bar{u}})+E(\mathbf{w}_{n})}{C(\mathbf{\bar
{u}})+C(\mathbf{w}_{n})}+\delta E(\mathbf{\bar{u}})+\delta E(\mathbf{w}%
_{n})\right] \\
&  \geq\lim\left[  \frac{E(\mathbf{\bar{u}})+E(\mathbf{w}_{n})}{\left\vert
C(\mathbf{\bar{u}})\right\vert +\left\vert C(\mathbf{w}_{n})\right\vert
}+\delta E(\mathbf{\bar{u}})+\delta E(\mathbf{w}_{n})\right] \\
&  \geq\lim\left[  \min\left(  \frac{E(\mathbf{\bar{u}})}{\left\vert
C(\mathbf{\bar{u}})\right\vert },\frac{E(\mathbf{w}_{n})}{\left\vert
C(\mathbf{w}_{n})\right\vert }\right)  +\delta E(\mathbf{\bar{u}})+\delta
E(\mathbf{w}_{n})\right]  .
\end{align*}
Now we consider two cases: first case $\frac{E(\mathbf{\bar{u}})}{\left\vert
C(\mathbf{\bar{u}})\right\vert }\geq\frac{E(\mathbf{w}_{n})}{\left\vert
C(\mathbf{w}_{n})\right\vert };$ then%
\begin{align*}
J_{\ast}  & \geq\lim\left[  \frac{E(\mathbf{w}_{n})}{\left\vert C(\mathbf{w}%
_{n})\right\vert }+\delta E(\mathbf{\bar{u}})+\delta E(\mathbf{w}_{n})\right]
\\
& =\lim\left[  J(\mathbf{w}_{n})+\delta E(\mathbf{\bar{u}})\right]  \geq
J_{\ast}+\delta E(\mathbf{\bar{u}}).
\end{align*}
This case cannot occur since it implies $\delta E(\mathbf{\bar{u}})\leq0$ and
this contradicts (\ref{ss}).

Then we have that
\[
\frac{E(\mathbf{\bar{u}})}{\left\vert C(\mathbf{\bar{u}})\right\vert }%
<\frac{E(\mathbf{w}_{n})}{\left\vert C(\mathbf{w}_{n})\right\vert }.
\]
In this case%
\begin{align*}
J_{\ast}  & \geq\lim\left[  \frac{E(\mathbf{\bar{u}})}{\left\vert
C(\mathbf{\bar{u}})\right\vert }+\delta E(\mathbf{\bar{u}})+\delta
E(\mathbf{w}_{n})\right] \\
& =\lim\left[  J(\mathbf{\bar{u}})+\delta E(\mathbf{w}_{n})\right]  \geq
J_{\ast}+\delta\lim E(\mathbf{w}_{n})
\end{align*}
Then%
\begin{equation}
\delta\lim E(\mathbf{w}_{n})\leq0.\label{ul}%
\end{equation}

Then by Lemma \ref{zeropiu} and (\ref{ul}) we have $\mathbf{w}_{n}%
\rightarrow0$ strongly in $X.$

$\square$

\textbf{Proof of Th. \ref{main-theorem}}. We shall use Theorem \ref{astra}.
Obviously assumptions (EC-1) and (EC-2) are satisfied with $G$ given by
(\ref{ggg}). Then by lemma \ref{GC} and Th. \ref{astra}, we have the existence
of soliton solutions. In order to prove that they form a family dependent of
$\delta$, it is sufficient to prove that $\delta_{1}\neq\delta_{2}$ in the
definition (\ref{J}) of $J$ implies $\mathbf{u}_{\delta_{1}}\neq
g\mathbf{u}_{\delta_{2}}$ for every $g\in G.$ We argue indirectly and assume
that $\mathbf{u}_{\delta_{1}}=g\mathbf{u}_{\delta_{2}}$ for some $g\in G.$
Then%
\[
\frac{E(g\mathbf{u}_{\delta_{2}})}{\left\vert C(g\mathbf{u}_{\delta_{2}%
})\right\vert }+\delta_{2}E(g\mathbf{u}_{\delta_{2}})=\frac{E(\mathbf{u}%
_{\delta_{1}})}{\left\vert C(\mathbf{u}_{\delta_{1}})\right\vert }+\delta
_{1}E(\mathbf{u}_{\delta1})
\]
and so, since $g\mathbf{u}_{\delta_{2}}=\mathbf{u}_{\delta_{1},}$
\begin{align*}
0 &  =\frac{E(g\mathbf{u}_{\delta_{2}})}{\left\vert C(g\mathbf{u}_{\delta_{2}%
})\right\vert }+\delta_{2}E(g\mathbf{u}_{\delta_{2}})-\left(  \frac
{E(\mathbf{u}_{\delta_{1}})}{\left\vert C(\mathbf{u}_{\delta_{1}})\right\vert
}+\delta_{1}E(\mathbf{u}_{\delta1})\right) \\
&  =\left(  \delta_{2}-\delta_{1}\right)  E(\mathbf{u}_{\delta_{1}}).
\end{align*}
Then, since $\delta_{1}\neq\delta_{2},$ $E(\mathbf{u}_{\delta_{1}})=0$ and so
$\mathbf{u}_{\delta_{1}}=0,$ which is a contradiction.

$\square$

\bigskip

\textbf{Proof of Th. \ref{imp}. }Since $\mathbf{u}_{\delta}=(u_{\delta
},v_{\delta})\in X=H^{2}(\mathbb{R})\times L^{2}(\mathbb{R})$ is a minimizer,
we have $J^{\prime}(\mathbf{u}_{\delta})=0.$Then%
\[
\frac{E^{\prime}(\mathbf{u}_{\delta})}{C(\mathbf{u}_{\delta})}-\frac
{E(\mathbf{u}_{\delta})}{C(\mathbf{u}_{\delta})^{2}}C^{\prime}(\mathbf{u}%
_{\delta})+\delta E^{\prime}(\mathbf{u}_{\delta})=0
\]
namely%
\[
\left(  C(\mathbf{u}_{\delta})+\delta C(\mathbf{u}_{\delta})^{2}\right)
E^{\prime}(\mathbf{u}_{\delta})=E(\mathbf{u}_{\delta})C^{\prime}%
(\mathbf{u}_{\delta}).
\]
Since, by (\ref{piu}), $C(\mathbf{u}_{\delta})>0,$ then $C(\mathbf{u}_{\delta
})+\delta C(\mathbf{u}_{\delta})^{2}>0,$ and hence we can divide both sides by
$C(\mathbf{u}_{\delta})+\delta C(\mathbf{u}_{\delta})^{2}$ and we get
\begin{equation}
E^{\prime}(\mathbf{u}_{\delta})=cC^{\prime}(\mathbf{u}_{\delta})\label{linda}%
\end{equation}
where%
\[
c=\frac{E(\mathbf{u}_{\delta})}{C(\mathbf{u}_{\delta})+\delta C(\mathbf{u}%
_{\delta})^{2}}.
\]
If we write (\ref{linda}) explicitely, we get for all $\varphi\in
H^{2}(\mathbb{R})$ and all $\psi\in L^{2}(\mathbb{R})$
\begin{align*}
\int\partial_{x}^{2}u_{\delta}\partial_{x}^{2}\varphi+W^{\prime}(u_{\delta
})\varphi & =c\int v_{\delta}\partial_{x}\varphi\\
\int v_{\delta}\psi & =c\int\psi\partial_{x}u_{\delta}%
\end{align*}
namely%
\begin{align*}
\partial_{x}^{4}u_{\delta}+W^{\prime}(u_{\delta})  & =-c\partial_{x}v_{\delta
}\\
v_{\delta}  & =c\partial_{x}u_{\delta}%
\end{align*}
and so we get%
\[
\partial_{x}^{4}u_{\delta}+c^{2}\partial_{x}^{2}u_{\delta}+W^{\prime
}(u_{\delta})=0
\]
Now. we can check directly that
\[
u(t,x)=u_{\delta}(x-ct)
\]
solves equation (\ref{1}) with initial conditions $\left(  u_{\delta
}(x),-c\partial_{x}u_{\delta}(x)\right)  .$

$\square$

\end{document}